# Asymptotic of summation functions

VICTOR VOLFSON

ABSTRACT We will study the asymptotic behavior of summation functions of a natural argument, including the asymptotic behavior of summation functions of a prime argument in the paper. A general formula is obtained for determining the asymptotic behavior of the sums of functions of a prime argument based on the asymptotic law of primes. We will show, that under certain conditions:

$$\sum_{p\le n} f(p) = \sum_{k=2}^{n} \frac{f(k)}{\log(k)}(1+o(1)),$$ where $p$ is a prime number.

In the paper, the necessary and sufficient conditions for the fulfillment of this formula are proved.

## 1. INTRODUCTION

An arithmetic function in the general case is a function $f$ defined on the set of natural numbers and taking values on the set of complex numbers. The name arithmetic function is due to the fact that this function expresses some arithmetic property of the natural series.

Summation is called a function of the form:

$$S(x) = \sum_{n\le x} f(n).$$

The Mertens function are the summation function - $M(\mathrm{x}) = \sum_{k\le x}\mu(k)$, where $\mu(k)$ is the Möbius arithmetic function.

Keywords: arithmetic function, summation function, Moebius function, Mertens function, asymptotic bound, Macleron's formula, Abel's summation formula, natural argument, prime argument, prime number

The Mobius function $\mu(\text{k}) = 1$, if a positive integer k has an even number of prime divisors of the first degree, $\mu(\text{k}) = -1$if a positive integer k has an odd number of prime divisors of the first degree and $\mu(\text{k}) = 0$ if a positive integer k has prime divisors of not only the first degree.

In addition to the summation functions of a natural argument, there are summation functions of a prime argument, where the summation is carried out only over primes. An example of the summation function of a prime argument is a function of the number of primes not exceeding the value $x$ - $\pi(x) = \sum_{p \le x} 1$.

An interesting problem is the study of the asymptotic behavior of summation functions.

We will study the asymptotic behavior of summation functions of a natural argument in Chapter 2 of the paper, and the asymptotic behavior of summation functions of a prime argument will be discussed in Chapter 3.

2. ASYMPTOTIC OF SUMMATION FUNCTIONS OF A NATURAL ARGUMENT

Let $S$ is a summation function of a natural argument. Then the relation is true:

$$S(n) = \sum_{k=1}^{n} f(k) = n \frac{\sum_{k=1}^{n} f(k)}{n} = nE[f,n] \ , \tag{2.1}$$

where $E[f,n]$is the average value (mathematical expectation) of the arithmetic function $f$ on the interval $[1,n)$.

Assertion 2.1

Let an arithmetic function $f$ take values: $a_1,...,a_l$ respectively with probabilities: $p_1,...,p_l (p_1 + ... + p_l = 1)$. Then it is executed:

$$S(n) = n \sum_{k=1}^{l} a_k p_k \ . \tag{2.2}$$

Proof

Based on the condition of the assertion:

$$E[f,n]=\sum_{k=1}^{l}a_k p_k\ . \qquad (2.3)$$

Having in mind (2.1) and (2.3) we get:

$$S(n)=n\sum_{k=1}^{l}a_k p_k\ ,$$

which corresponds to (2.2).

As an example, we consider the Mertens function.

Given the above, the values of the Möbius arithmetic function are: $a_1=1, a_2=-1, a_3=0$. Having in mind [1], the probabilities are respectively equal to: $p_1=p_2=3/\pi^2+\mathrm{o}(1), p_3=1-6/\pi^2+o(1)$. Therefore, based on assertion 2.1:

$$M(n)=\mathrm{n}(3/\pi^2+o(1)-3/\pi^2+o(1))=o(n)\,.$$

Now we consider the summation functions of quantity $Q(n)=\sum_{n\in A}1$, which determine the number of members of the natural series that satisfy the condition $A$.

Thus, in this case, the arithmetic function takes only two values: $f(n)=1$ with the probability of fulfilling the condition $A$ - $p(n)$ and $f(n)=0$ with probability $1-p(n)$. Therefore, for the summation function of the quantity, based on assertion 2.1, the following statement holds.

Assertion 2.2

It is performed for the summation function of the quantity:

$$Q(n)=np(n)\,. \qquad (2.4)$$

As an example, we consider the asymptotic upper bound for the summation function of the number of primes not exceeding the value $n$ - $\pi(n)$.

It was proved in [2] that the probability of a natural number from the interval $[2,n]$ to be prime is equal to:

$$p(n)=(1+o(1))/\log(n)\,. \qquad (2.5)$$

Having in mind (2.4) and (2.5) we obtain:

$$\pi(\mathrm{n}) = n(1+o(1)) / \log(n) . \qquad (2.6)$$

It is possible to determine the asymptotic behavior of the summation function using the Euler – Maclaurin formula [3] in the case when the function $f$ is elementary and sufficiently smooth (has the required number of derivatives) on the interval $[1,n]$.

In particular, the asymptotic behavior of the summation function is determined by the formula (if the function $f$ on the interval $[1,n]$ is strictly decreasing and $\lim_{n\to\infty} f(n) = 0$ ):

$$S(n) = \sum_{k=1}^{n} f(k) = \int_{t=1}^{n} f(t)dt + C + O(f(n)) , \qquad (2.7)$$

where $C$ is the constant.

As an example of using (2.7), we consider the definition of the asymptotic of the following summation function:

$$S(n) = \sum_{k=1}^{n} \frac{\log(k)}{k} = \int_{t=1}^{n} \frac{\log(t)\,\mathrm{dt}}{t} + C + O(\frac{\log(n)}{n}) = \log\log(n) + C + O(\frac{\log(n)}{n}) . \qquad (2.8)$$

If the function $f$ on the interval $[1,n]$ is non-decreasing, then the asymptotic behavior of the summation function is determined by the formula:

$$S(n) = \sum_{k=1}^{n} f(k) = \int_{t=1}^{n} f(t)dt + O(f(n)) . \qquad (2.9)$$

As an example of using (2.9), we consider the definition of the asymptotics of the following adder function:

$$S(n) = \sum_{k=1}^{n} k^m = \int_{t=1}^{n} t^m dt + O(n^m) = \frac{n^{m+1}}{m+1} + O(n^m) . \qquad (2.10)$$

3. ASYMPTOTIC OF SUMMATION FUNCTIONS OF A PRIME ARGYMENT

The summation function of a prime argument is an arithmetic function of the form:

$$S(n) = \sum_{p \le n} f(p) , \qquad (3.1)$$

where $p$ is a prime number.

Conjecture 3.1

If the function $f(k)/\log(k)$ is smooth enough, i.e. has the desired number of derivatives on the interval $[2,n]$, then:

$$S(\mathrm{n}) = \sum_{p\le n} f(p) = \sum_{k=2}^{n} \frac{f(k)}{\log(k)}(1+o(1))$$

and the asymptotic $\sum_{p\le n} f(p)$ can be determined using the Euler-Maclaurin formula.

Justification

Note that a sequence of values $f(p)$: $f(2), f(3), f(5),\ldots$ is obtained by sifting a sequence of values $f(n)$: $f(1), f(2), f(3), f(4), f(5),\ldots$.

Thus, the arithmetic function $f(p) = f(n)$, if $n = p$, and $f(p) = 0$, if $n \ne p$.

Earlier, we showed that the probability of a natural number to be a prime in the interval $[2,n]$ is determined by the formula (2.5). We replace the function $f(p)$ with its mathematical expectation on the interval $[2,n]$:

$$S(\mathrm{n}) = \sum_{p\le n} f(p) = \sum_{k=2}^{n} \frac{f(k)}{\log(k)}(1+o(1)). \tag{3.2}$$

Since the function $f(k)/\log(k)$ is quite smooth, i.e. has the desired number of derivatives on the interval $[2,n]$, then the Euler-Maclaurin formula can be used to determine the asymptotic $\sum_{p\le n} f(p)$.

Let's look at examples of determining the asymptotic of summation functions of a prime argument using (3.2).

Having in mind (2.7) and (3.2) we get:

$$\sum_{p\le n} \frac{1}{p} = \sum_{k=2}^{n} \frac{1}{k\log(k)}(1+o(1)) = \log\log(n)(1+o(1)). \tag{3.3}$$

Based on (2.9) and (3.2) we obtain:

$$\sum_{p\le n}\log(p)=\sum_{k=2}^{n}\frac{\log(k)}{\log(k)}(1+o(1))=n(1+o(1)). \quad (3.4)$$

Having in mind (2.7) and (3.2) we get:

$$\sum_{p\le n}\frac{\log(p)}{p}=\sum_{k=2}^{n}\frac{\log(k)}{k\log(k)}(1+o(1))=\log(n)(1+o(1)). \quad (3.5)$$

The results (3.3), (3.4) and (3.5) correspond to [3], [4].

Based on (3.2), one can also obtain other asymptotic estimates.

Having in mind (2.9) and (3.2) we get:

$$\sum_{p\le n}p^{l}\log(\mathrm{p})\approx\sum_{k=2}^{n}\frac{k^{l}\log(k)}{\log(k)}=\frac{n^{l+1}}{l+1}(1+o(1)), \quad (3.6)$$

where $l\ge 0$.

Based on (2.7) and (3.2) we obtain:

$$\sum_{p\le n}\frac{\log^{l}p}{p}\approx\sum_{k=2}^{n}\frac{\log^{l}k}{k\log(k)}(1+o(1))=\frac{\log^{l}n}{l}(1+o(1)), \quad (3.7)$$

where $l\ge 1$.

Now consider the requirements under which conjecture 3.1 is true.

Let $a_k=1$, if $k$ is prime and $a_k=0$, if $k$ compound and let $b_k=\frac{1}{\ln k}, b_1=0$.

Denote: $A(\mathrm{n})=\sum_{k\le n}a_k$ and $B(n)=\sum_{k\le n}b_k$.

It is required that:

$$\lim_{n\to\infty}\frac{A(n)}{B(n)}=1. \quad (3.8)$$

It follows from (3.8) that in order to fulfill conjecture 3.1 it is necessary:

$$\lim_{n\to\infty}\frac{\sum_{k=1}^{n} a_k f(k)}{\sum_{k=1}^{n} b_k f(k)} = 1. \qquad (3.9)$$

Assertion 3.2

If the conditions are met:

1. $\lim_{n\to\infty}\dfrac{\int_1^n B(t) f'(t)dt}{B(n) f(n)}$ is not equal to 1.

2. $f(x)$ is monotonous and $f'(x) \neq 0$.

3. $\lim_{n\to\infty}\int_1^n B(t) f'(t)dt = \pm\infty$,

then:

$$\lim_{n\to\infty}\frac{\sum_{k=1}^{n} a_k f(k)}{\sum_{k=1}^{n} b_k f(k)} = 1.$$

Proof

Based on Abel's summation formula:

$$\sum_{k=1}^{n} a_k f(k) = A(n) f(n) - \int_1^n A(t) f'(t)dt, \sum_{k=1}^{n} b_k f(k) = B(n) f(n) - \int_1^n B(t) f'(t)dt. \qquad (3.10)$$

Having in mind (3.9) and (3.10) we get:

$$\frac{\sum_{k=1}^{n} a_k f(k)}{\sum_{k=1}^{n} b_k f(k)} = \frac{A(n) f(n) - \int_1^n A(t) f'(t)dt}{B(n) f(n) - \int_1^n B(t) f'(t)dt} = \frac{A(n)}{B(n)} \cdot \frac{1 - \dfrac{\int_1^n A(t) f'(t)dt}{A(n) f(n)}}{1 - \dfrac{\int_1^n B(t) f'(t)dt}{B(n) f(n)}}. \qquad (3.11)$$

Let the conditions be satisfied:

1. $\lim_{n\to\infty}\frac{\int_1^n B(t)f'(t)dt}{B(n)f(n)}$ is not equal to 1.

2. Let $f(x)$ is monotonous and $f'(x)\neq 0$.

3. $\lim_{n\to\infty}\int_1^n B(t)f'(t)dt=\pm\infty$.

We show that when these conditions are met:

$$\lim_{n\to\infty}\frac{\int_1^n A(t)f'(t)dt}{A(n)f(n)}=\lim_{n\to\infty}\frac{\int_1^n B(t)f'(t)dt}{B(n)f(n)}. \tag{3.12}$$

which corresponds to (3.12).

Indeed, using the Local rule, we get:

$$\lim_{n\to\infty}\frac{\dfrac{\int_1^n A(t)f'(t)dt}{A(n)f(n)}}{\dfrac{\int_1^n B(t)f'(t)dt}{B(n)f(n)}}=\lim_{n\to\infty}\frac{B(n)}{A(n)}\cdot\lim_{n\to\infty}\frac{\int_1^n A(t)f'(t)dt}{\int_1^n B(t)f'(t)dt}=\lim_{n\to\infty}\frac{A(n)f'(n)}{B(n)f'(n)}=1, \tag{3.13}$$

which corresponds to (3.12).

Then, based on (3.8) and (3.13), we obtain:

$$\lim_{n\to\infty}\frac{\sum_{k=1}^n a_k f(k)}{\sum_{k=1}^n b_k f(k)}=\lim_{n\to\infty}\frac{A(n)}{B(n)}\cdot\frac{1-\dfrac{\int_1^n A(t)f'(t)dt}{A(n)f(n)}}{1-\dfrac{\int_1^n B(t)f'(t)dt}{B(n)f(n)}}=1,$$

which is consistent with the assertion 3.2.

Corollary 3.3

Conditions (1), (3) in assertion 2.2 correspond to:

1. $\lim_{n\to\infty}\dfrac{\int_2^n \frac{tf'(t)}{\log(t)}dt}{\frac{nf(n)}{\log(n)}}$ is not equal to 1.

3. $\lim_{n\to\infty}\int_2^n \frac{tf'(t)}{\log(t)}dt = \pm\infty$.

Proof

Based on the asymptotic law of primes:

$$\mathrm{A}(n) = B(n) = \frac{n}{\log(n)}(1+o(1)). \tag{3.14}$$

Substituting (3.14) in conditions 1 and 3 we get:

1. $\lim_{n\to\infty}\dfrac{\int_2^n \frac{tf'(t)}{\log(t)}dt}{\frac{nf(n)}{\log(n)}} \neq 1$.

3. $\lim_{n\to\infty}\int_2^n \frac{tf'(t)}{\log(t)}dt = \pm\infty$.

Let's look at an example: $\sum_{p\le n}\frac{1}{p}$. In this case $f(n) = \frac{1}{n}$.

1. $\lim_{n\to\infty}\dfrac{-\int_2^n \frac{dt}{\mathrm{t}^3\log(t)}}{\frac{n}{\mathrm{n}\log(n)}} = -\lim_{n\to\infty}\frac{n\log(n)}{n^4\log(n)} = 0 \neq 1$.

2. $f(n) = \frac{1}{n}$ is monotonic function and $f'(n) = -\frac{1}{n^2} \neq 0$.

3. $-\lim_{n\to\infty}\int_2^n \frac{tdt}{\mathrm{t}^2\log(t)} = -\lim_{n\to\infty}(\log\log(n)) = -\infty$.

Thus, all the conditions of Assertion 3.2 and Corollary 3.3 are satisfied in this case; therefore, asymptotic formula (3.3) is true.

Similarly to asymptotic formula (3.3), we can verify that the conditions of Assertion 3.2 and Corollary (3.3) are satisfied for the asymptotic formulas (3.4), (3.5), (3.6) and (3.7).

The conditions of Assertion 3.2 and Corollary 3.3 are sufficient to satisfy Conjecture 3.1.

Let us consider a special case of Assertion 3.2 and Corollary 3.3, when the function $f$ is monotonically increasing and tends to infinity.

Assertion 3.4

Let $\lim_{n\to\infty} f(n) = \infty$, $f'(n)$ is a continuous function, $f'(n) > 0$ and $\lim_{n\to\infty} \frac{f(n)}{nf'(n)} \neq 0$. Then all the conditions of Assertion 3.2 and Corollary 3.3 are satisfied.

Proof

Check the 1st condition:

$$\lim_{n\to\infty} \frac{\int_2^n \frac{tf'(t)}{\log(t)} dt}{\frac{nf(n)}{\log(n)}} = \lim_{n\to\infty} \frac{\frac{nf'(n)}{\log(n)}}{(\frac{nf(n)}{\log(n)})'} = \lim_{n\to\infty} \frac{\frac{nf'(n)}{\log(n)}}{\frac{nf'(n)}{\log(n)} + (\frac{n}{\log(n)})' f(n)} = \lim_{n\to\infty} \frac{1}{1 + \frac{f(n)}{nf'(n)}} \neq 1$$

Check the 2nd condition:

$f(n)$ is monotonic function and $f(n) \neq 0$.

Check the 3rd condition:

$\int_2^n \frac{tf'(t)}{\log(t)} dt \geq \int_2^n \frac{2f'(t)}{\log(2)} dt = \frac{2}{\log(2)} \int_2^n f'(t) dt = \frac{2}{\log(2)} (f(n) - f(2))$ function increases unlimitedly, like $f(n)$.

Assertion 3.4 holds for the asymptotic formulas (3.4), (3.6), and also for the function $\sum_{p\leq n} p^m$ where $m > -1$.

Assertion 3.5

Suppose that the function $f'$ is continuous and the above conditions are satisfied, then the following asymptotic equality holds:

$$\int_2^n \frac{f(t)dt}{\log(t)} = (\frac{f(n)n}{\log(n)} - \frac{2f(2)}{\log(2)} - \int_2^n \frac{tf'(t)dt}{\log(t)})(1+o(1)) . \qquad (3.15)$$

Proof

We use that: $J_1(t) = \int_2^t \frac{du}{\log(u)} = \frac{t}{\log(t)}(1+o(1))$ and integrating by parts we get:

$$\int_2^n \frac{f(t)dt}{\log(t)} = f(t)J_1(t)\,|_2^n - \int_2^n J_1(t) f'(t)dt = (\frac{f(n)n}{\log(n)} - \frac{2f(2)}{\log(2)} - \int_2^n \frac{tf'(t)dt}{\log(t)})(1+o(1)) ,$$

which corresponds to (3.15).

Let's look at an example using assertion 3.5.

We prove the asymptotic equality:

$$\sum_{p \le n} p^m = \frac{n^{m+1}}{(m+1)\log(n)}(1+o(1)) , \qquad (3.16)$$

with values $m > -1$.

Based on assertion 3.5:

$$J = \int_2^n \frac{t^m dt}{\log(t)} = \frac{n^{m+1}}{\log(n)} - \frac{2^{m+1}}{\log(2)} - mJ . \qquad (3.17)$$

Having in mind (3.17):

$$J = \int_2^n \frac{t^m dt}{\log(t)} = \frac{1}{m+1}(\frac{n^{m+1}}{\log(n)} - \frac{2^{m+1}}{\log(2)}) . \qquad (3.18)$$

Based on (3.18) we get:

$$\sum_{p \le n} p^m = \sum_{k=2}^n \frac{k^m}{\log(k)}(1+o(1)) = J(1+o(1)) = \frac{n^{m+1}}{(m+1)\log(n)}(1+o(1)) ,$$

which corresponds to (3.16).

The asymptotic equality (3.16) corresponds to [3].

We now consider the necessary condition for fulfilling hypothesis 3.1.

We use the notation of assertion 3.2.

Assertion 3.6

Let:

$$\lim_{n\to\infty}\frac{\sum_{k=1}^{n}a_k f(k)}{\sum_{k=1}^{n}b_k f(k)}=1$$

Then, for $p\to\infty$ ( $p$ is a prime number), the following holds:

$$\left|\frac{f(p)}{\sum_{k=1}^{p}b_k f(k)}\right|\to 0\,. \qquad (3.19)$$

Proof

If:

$$\lim_{n\to\infty}\frac{\sum_{k=1}^{n}a_k f(k)}{\sum_{k=1}^{n}b_k f(k)}=1$$

Then holds:

$$\lim_{n\to\infty}\left|\frac{\sum_{k=1}^{n}a_k f(k)}{\sum_{k=1}^{n}b_k f(k)}-\frac{\sum_{k=1}^{n-1}a_k f(k)}{\sum_{k=1}^{n-1}b_k f(k)}\right|=0\,. \qquad (3.20)$$

Take $n=p$, where $p$ is a prime. Then, based on (3.20):

$$\lim_{p\to\infty}\left|\frac{\sum_{k=1}^{p}a_k f(k)}{\sum_{k=1}^{p}b_k f(k)}-\frac{\sum_{k=1}^{p-1}a_k f(k)}{\sum_{k=1}^{p-1}b_k f(k)}\right|=0 \qquad (3.21)$$

Let us denote (at $k=p$) the value $a_k f(k)=a_p f(p)$.

We transform (3.21) and taking into account that $a_p=1$ we get:

$$\left|\frac{\sum_{k=1}^{p} a_k f(k)}{\sum_{k=1}^{p} b_k f(k)} - \frac{\sum_{k=1}^{p-1} a_k f(k)}{\sum_{k=1}^{p-1} b_k f(k)}\right| = \left|\frac{a_p f(p)\sum_{k=1}^{p-1} b_k f(k) - b_p f(p)\sum_{k=1}^{p-1} a_k f(k)}{\sum_{k=1}^{p} b_k f(k)\sum_{k=1}^{p-1} b_k f(k)}\right| = | f(p) | \left|\frac{1 - b_p \frac{\sum_{k=1}^{p-1} a_k f(k)}{\sum_{k=1}^{p-1} b_k f(k)}}{\sum_{k=1}^{p} b_k f(k)}\right|$$

Since when $p \to \infty$ the value $b_p \to 0$ and $\lim_{n\to\infty} \frac{\sum_{k=1}^{n} a_k f(k)}{\sum_{k=1}^{n} b_k f(k)} = 1$, it is executed:

$$\left|\frac{f(p)}{\sum_{k=1}^{p} b_k f(k)}\right| \to 0 ,$$

which corresponds to (3.19).

For example, the specified necessary condition is not satisfied for the function $f(p) = 2^p$, since $\sum_{k=1}^{p} b_k f(k) = \sum_{k=1}^{p} \frac{2^k}{\log(k)} = \frac{2^{p+1}}{\log p}$.

Now we consider the asymptotic behavior of the mean value (mathematical expectation) of the arithmetic function of a prime argument.

Corollary 3.7

Based on formula (2.1), under the conditions of satisfying the sufficient conditions of Assertion 3.2, it is executed for the average value (mathematical expectation) of the arithmetic function of a prime argument:

$$E[f(p), n] = \sum_{p \le n} f(p) / n = \frac{1}{n}\sum_{k=2}^{n} \frac{f(k)}{\log(k)}(1 + o(1)) . \qquad (3.22)$$

As an example, using (3.4) and (3.22), we determine the average value for the following arithmetic function of a prime argument:

$$E[\log(p), n] = \sum_{p \le n} \log(p) / \mathrm{n} = \frac{1}{n}\sum_{k=2}^{n} \frac{\log(k)}{\log(k)}(1 + o(1)) = 1 + o(1) .$$

## 4. CONCLUSION AND SUGGESTIONS FOR FURTHER WORK

The next article will continue to study the behavior of some sums.

## 5. ACKNOWLEDGEMENTS

Thanks to everyone who has contributed to the discussion of this paper. I am grateful to everyone who expressed their suggestions and comments in the course of this work.

References

1. Volfson V.L. Investigation of the asymptotic behavior of the Mertens function, Applied Physics and Mathematics No. 6, 2017, pp. 46-50.

2. Volfson V.L. Asymptotics of the greatest distance between adjacent primes and the Hardy-Littlewood conjecture, Applied Physics and Mathematics No. 2, 2020, pp. 38-44.

3. K. Prahar, Distribution of primes, M, Mir, 1967 -512 p.

4. A.A. Buchstaff. "Theory of numbers", From the "Enlightenment", Moscow, 1966, 384 p.